\documentclass[conference]{IEEEtran}
\IEEEoverridecommandlockouts
\usepackage{cite}
\usepackage{amsmath,amssymb,amsfonts}
\usepackage{algorithmic}
\usepackage{graphicx}
\usepackage{textcomp}
\usepackage{microtype}
\usepackage{xcolor}
\def\BibTeX{{\rm B\kern-.05em{\sc i\kern-.025em b}\kern-.08em
    T\kern-.1667em\lower.7ex\hbox{E}\kern-.125emX}}

\usepackage{amsmath}
\DeclareMathOperator{\trace}{Tr}

\usepackage{tikz}
\usepackage{textcomp}
\usepackage{hyperref}
\usepackage{lipsum}

\newcommand\copyrighttext{%
  \footnotesize \textcopyright 2018 IEEE. Personal use of this material is permitted.
  Permission from IEEE must be obtained for all other uses, in any current or future
  media, including reprinting/republishing this material for advertising or promotional
  purposes, creating new collective works, for resale or redistribution to servers or
  lists, or reuse of any copyrighted component of this work in other works.
  DOI: \href{https://dx.doi.org/10.1109/ScalA.2018.00014}{10.1109/ScalA.2018.00014}}
\newcommand\copyrightnotice{%
\begin{tikzpicture}[remember picture,overlay]
\node[anchor=south,yshift=10pt] at (current page.south) {\fbox{\parbox{\dimexpr\textwidth-\fboxsep-\fboxrule\relax}{\copyrighttext}}};
\end{tikzpicture}%
}

\begin{document}

\title{On Advanced Monte Carlo Methods for Linear Algebra on Advanced Accelerator Architectures\\
}

\author{
\IEEEauthorblockN{1\textsuperscript{st} Anton Lebedev}
\IEEEauthorblockA{\normalsize Institute for Theoretical Physics,\\
University of T\"{u}bingen,\\
Germany\\
Email: anton.lebedev@student.\\
uni-tuebingen.de}
\and
\IEEEauthorblockN{2\textsuperscript{nd} Vassil Alexandrov}
\IEEEauthorblockA{\normalsize ICREA, Catalan Institution for\\
Research and Advanced Studies\\
Barcelona Supercomputing Centre, Spain\\
Email: vassil.alexandrov@bsc.es}
}
\maketitle

\copyrightnotice

\begin{abstract}
In this paper we present  computational experiments with the Markov Chain Monte Carlo Matrix Inversion ($(\text{MC})^2\text{MI}$)
on several accelerator architectures  and investigate their impact on performance and scalability of the method.
The method is used as a preconditioner and for solving the corresponding system of linear equations  iterative methods, such as generalized minimal residuals (GMRES) or bi-conjugate gradient (stabilized) (BICGstab), are used.

Numerical experiments are carried out to highlight the benefits and deficiencies of both approaches and to assess their overall
usefulness in light of scalability of the method.
\end{abstract}


\section{Introduction}

Solving systems of linear algebraic equations (SLAE) in the form of $Bx = b$ or inverting a matrix $B$ is of unquestionable importance in many scientific fields. Iterative solvers are used widely to compute the solutions of these systems and such approaches are often the method of choice due to their predictability and reliability when considering accuracy and speed. They, however, may become prohibitive for large-scale problems as they can be very time consuming to compute.  The complexity of these methods, in the serial case,  is $O(kn^2)$ for dense matrices in the iterative methods case and $O(n^3)$ for direct methods with dense matrices while solving SLAE if common elimination or annihilation schemes (e.g. Gaussian elimination, Gauss-Jordan methods) are employed \cite{golub1996matrix}. Therefore, these algorithms often rely on preconditioners to speed up the computations and/or to ensure faster convergence.

Monte Carlo (MC) methods complexity is linear in matrix size \cite{AlexandrovStrassburgA14}, \cite{AlexandrovV15} and can quickly yield a rough estimate of the solution
by sampling a random variable whose mathematical expectation is the desired solution. For some problems an estimate is sufficient or even favourable, due to the accuracy of the underlying data. Therefore, it should be pointed out, that  Monte Carlo methods may  be efficiently used as preconditioners. 

Depending on the method used to compute the preconditioner, the savings and end-results vary. A very sparse preconditioner may be computed quickly, but it is unlikely to
greatly reduce the run time to solution. On the other hand, computing a rather dense preconditioner is computationally expensive and might be time or cost prohibitive. Therefore, finding a good preconditioner that is computationally efficient, while still providing substantial improvement to the iterative solution process, is a worthwhile research topic.

A variety of parallel Monte Carlo methods have been developed within the past 20 years. A comprehensive compendium of the Monte Carlo functions and strategies of parallelization can be found in \cite{vajargah2001parallel, strassburgPhd2014, branford2008, AlexandrovStrassburgA14, AlexandrovV15} .

In this work we present an enhanced version of a SPAI (SParse Approximate Inverse) preconditioner that is based on parallel  Monte Carlo methods presented in  \cite{AlexandrovStrassburgA14} and \cite{AlexandrovV15}. This new optimized version is compared against the previous one, taken as a baseline, as well as against 
MSPAI, which is the main accepted deterministic algorithm for SPAI preconditioning. Our results show that the Monte Carlo-based algorithm can be used instead of MSPAI to reduce the computation time and resource usage while producing results with  similar or better quality.

Also a scalability analysis is carried out, showing that the random patterns in the memory access have a strong influence on the performance of the algorithm. Further research, to solve this issues, is proposed within the context of quasi-Monte Carlo Methods.

The next section gives and overview of related work. Monte Carlo methods, and the specific matrix inversion algorithm that is discussed as a SPAI preconditioner, are presented in section \ref{section_algorithm}. Section~\ref{section:paralleldetails}  presents parallel  approach of the MonteCarlo and the hybrid algorithm. Section \ref{section:paralleldetails} shows the approach  and  methodology applied in the enhancement of the parallel implementations
Sections \ref{sec:AlgModifications} and \ref{sec:NumExp} present corresponding results and analysis of the implementations. The conclusion \ref{section_conclusion}  summarises the results and  outlines the future work.

\section{Related Work}
\label{section_related_work}

Research efforts in the past have been directed towards optimizing the approach of
sparse approximate inverse preconditioners. Improvements to the computation of the Frobenius norm have been
proposed for example by concentrating on sparse pattern selection
strategies~\cite{carpentieri2000some}, or building a symmetric preconditioner by
averaging off-diagonal entries~\cite{carpentieri2000experiments}. Further, it has been
shown that the sparse approximate inverse preconditioning approach is also a viable
course of action on large-scale dense linear systems~\cite{alleon1997sparse}. This is of
special interest to us, as the Monte Carlo method we are proposing in this paper is part
of a bigger family. It includes serial and parallel Monte Carlo algorithms for the
inversion of sparse, as well as dense matrices, and the solution of systems of linear
algebraic equations. The proposed Monte Carlo algorithm has been developed and enhanced
in the last decades, and several key advances in serial and parallel Monte Carlo
methods for solving such problems have been made~\cite{ Branford2008grid, Dimov1998convergent, Dimov2007}.
There has been an increased research interest in parallel Monte Carlo
methods for linear algebra in the past few years, and recent example is the Monte Carlo
Synthetic Acceleration (MCSA) developed through MCREX project at ORNL\cite{mcrex2013}.
Future work that deals with a parallel implementation of the presented algorithm is
being considered further in this section and in section~\ref{section_algorithm}.

In the past there have been differing approaches and advances towards a parallelisation
of the SPAI preconditioner. In recent years the
class of Frobenius norm minimizations that has been used in the original SPAI
implementation~\cite{benzi1996sparse} was modified and is provided in a parallel SPAI
software package. One implementation of it, by the original authors of SPAI, is the
Modified SParse Approximate Inverse (MSPAI~\cite{Huckle2010}).

This version provides a class of modified preconditioners such as MILU (modified ILU),
interface probing techniques and probing constraints to the original SPAI, apart from a
more efficient, parallel Frobenius norm minimization. Further, this package also
provides two novel optimization techniques. One option is using a dictionary in order to
avoid redundant calculations, and to serve as a lookup table. The second option is 
an option  to switch to a less computationalyl intensive, sparse QR
decomposition whenever possible. This optimized code runs in parallel, together with a
dynamic load balancing.

\subsection{Using SParse Approximate Inverse as Preconditioner (SPAI)}

The SPAI algorithm~\cite{grote2006spai} is used to compute a sparse approximate inverse
matrix $M$ for a given sparse input matrix $B$. This is done by minimizing $\Vert BM-I\Vert_F$. 
The algorithm explicitly computes the approximate inverse, which is
intended to be applied as a preconditioner of an iterative method. The SPAI application
provides the option to fix the sparsity pattern of the approximate inverse a priori or
capture it automatically. Since the introduction of the original SPAI in 1996, several advances, building upon the initial implementation, have been made. Two newer implementations are provided by the original authors, the aforementioned MSPAI, and the highly scalable Factorized SParse Approximate Inverse (FSPAI~\cite{huckle2003factorized}). The intended use of both differs depending on the problem at hand. Whereas MSPAI is used as a preconditioner for large sparse and ill-conditioned systems of linear equations, FSPAI is applicable only to symmetric positive definite systems of this kind. FSPAI is based around an inherently parallel implementation, generating the approximate inverse of the Cholesky factorization for the input matrix. MSPAI on the other hand is using an extension of the well-known Frobenius norm minimization that has been introduced in the original SPAI. 

The algorithm attempts to solve a system of linear equations of the form $Bx=b$. Its input is a sparse, square coefficient matrix $B$. The right hand side vector $b$ can either be provided by the user, or is arbitrarily defined by the software implementation. In the case of the SPAI application suite, if no right hand side vector is handed to the algorithm, it constructs one by multiplying matrix $B$ with a vector consisting of all ones. In a general case, an input matrix $B$ is passed to SPAI as a file. The program then computes a preconditioner using the Frobenius norm, afterwards it uses this intermediate result as an input to an appropriate solver.


\section{Monte~Carlo Approach}
\label{section_algorithm}

Monte~Carlo methods are probabilistic methods that use random numbers to either
simulate a stochastic behaviour or to estimate the solution of a problem. They are good
candidates for parallelisation due to the fact that, in principle, many independent samples are
used to estimate the solution. These samples can be calculated in parallel, thereby
speeding up the solution finding process. The so  designed  and developed parallel Monte~Carlo
methods possess the following main generic properties  \cite{AlexandrovStrassburgA14, AlexandrovV15}: efficient distribution of the compute data, minimum communication during the computation
and increased precision being achieved by adding extra refinement computations. Consideration of all these properties naturally leads to scalable algorithms.
It has to be noted that the quality of the solutions obtained using a Monte~Carlo method is dependent upon the availability
of independent (pseudo) random numbers, which is a concern in parallel environments.

\subsection{Algorithm}\label{subsec:Algorithm}
The following procedure has been presented in~\cite{branford2008}  and allows to extend
the Monte Carlo algorithm for processing diagonally dominant matrices, that is used as the foundation for this work (c.f.~\cite{strassburg2013scalability}),  to the case of general matrices  \cite{AlexandrovStrassburgA14} \cite{AlexandrovV15}.

Let us recall for simplicity the key details from \cite{strassburg2013scalability, AlexandrovStrassburgA14, AlexandrovV15}. We assume the general case where $\Vert B\Vert > 1$, with $\Vert\cdot\Vert$ being an arbitrary matrix norm, and consider the splitting
\begin{equation}\label{f-mcla-mca-mcm-0}
B = \hat{B} - C,
\end{equation}
where the off-diagonal elements of $\hat{B}$ are the same as those of $B$, and the diagonal elements of $\hat{B}$ are defined as $\hat{b}_{ii}=b_{ii}+\alpha_{i} \Vert B\Vert$, choosing in most cases $\alpha_{i}> 1$ for $i=1,2,...,n$. For the simplicity of the algorithm it is often easier to fix single  $\alpha$ . 
In the general case, $\Vert B\Vert > 1$, make the initial split $B = \hat{B} - C$. From this compute $A = B_1^{-1}B_2$, $B_1 = diag(\hat{B})$ which satisfies $\Vert A\Vert < 1$.
Then the inverse of $\hat{B}$  is generated by 
\begin{equation}
[\hat{B}^{-1}]_{rr^{\prime}}\approx \frac{1}{N} 
\sum_{s=1}^N\left[\sum_{( j| s_{j}=r^{\prime})} 
W_{j}\right],\label{f-mcla-mca-dd-12}
\end{equation}
where $( j | s_{j}=r^{\prime})$ means that only
\[ W_{j} = \frac{a_{r s_1} a_{s_1 s_2} \dots a_{s_{j-1} s_{j}}}{p_{r s_1} p_{s_1 s_2} \dots p_{s_{j-1} s_{j}}}, \]
for which $s_{j} = r^{\prime}$ are included in the sum (\ref{f-mcla-mca-dd-12}).
Calculating $\Vert B\Vert$ can be an expensive operation, so any \textit{a-priori} information
allowing for a reasonable estimate here is useful but not strictly necessary.
From this it is then necessary to work back and recover $B^{-1}$ from $\hat{B}^{-1}$. To
do this an iterative process ($k = n-1, n-2,\ldots , 0$) is used on $\hat{B}^{-1}$:
\begin{equation}\label{f-mcla-mca-mcm-1} B^{-1}_{k} = B^{-1}_{k+1} + \frac{B^{-1}_{k+1}S_{k+1}B^{-1}_{k+1}}{1 - \trace \left(B^{-1}_{k+1}S_{k+1}\right)}, \end{equation}
where $B^{-1}_n = \hat{B}^{-1}$ and $S_{i}$ is all zero except for the $\{ii\}^{th}$
component, which is from the matrix $S=\hat{B} - B$. Then $B^{-1}=B^{-1}_0 $.

The make up of matrix $S$ means that while (\ref{f-mcla-mca-mcm-1}) looks complicated it
is, in fact simply an update of the matrix by a scaled outer product of the $(k+1)^{th}$ column with the $(k+1)^{th}$ row. 
There are obvious simplifications possible to ensure that
many multiplications by zero are not performed. This method of splitting and recovery
leads to Algorithm~1 \cite{AlexandrovStrassburgA14}, which details a MC algorithm for inverting
general matrices and is given below for completeness. Further details on the recovery of the original inverse can be found in \cite{Fathi07}.
\\

\noindent \textbf{Algorithm 1: } Monte Carlo Algorithm for Inverting General Matrices
  \begin{enumerate} 
    \item Read in matrix $B$
      \begin{enumerate} 
	\item Input matrix $B$, parameters $\varepsilon$ and $\delta$
      \end{enumerate}
    \item Remove a set percentage of the smallest (in magnitude) entries of the matrix.
    \item Calculate intermediate matrices ($\hat{B}$, $B_1$)
      \begin{enumerate} 
	\item Split $B = \hat{B} - (\hat{B} - B)$, where $\hat{B}$ is a diagonally
	dominant matrix
      \end{enumerate}
    \item Apply the algorithm for inverting diagonally dominant
   matrices from~\cite{strassburg2013scalability} with $B=\hat{B}$ to obtain
   $\hat{B}^{-1}$
   \item Recovery of $B^{-1}$ from $\hat{B}^{-1}$
    \begin{enumerate} 
     \item Compute $S = \hat{B} - B$
      \begin{enumerate} 
       \item Let $S_i$ for $i = 1, 2, \ldots, n$ where each $S_i$ has just
        one of the n[on-zero elements of the matrix $S$
       \item Set $B_n^{-1} = \hat{B}^{-1}$
       \item  Apply $B_{i-1}^{-1} = B_i^{-1} + \frac{B_i^{-1}S_{i}B_i^{-1}}{1 - \trace\left(B_i^{-1}S_i\right)}$ for $i = n, n-1, \ldots, 1$
     \end{enumerate}
    \item Then $B^{-1} = B_{0}^{-1}$
   \end{enumerate}
  \end{enumerate}

Note that the second step is optional and is relevant only when a reduction in the amount of data being communicated is desired. Its
influence has been investigated and the results are presented in sec. \ref{subsec:FitnessForPurpose}.

The above algorithm was modified to develop an MPI version of the algorithm. Several enhancements of the algorithm, as well as modifications concerning GPU implementation, are listed in the next section and were able to substantially improve its performance in generating rough inverses of the input matrices. The result can then be used directly as a preconditioner for solving a system of linear algebraic equations or further improved. We propose the use of an iterative refinement process, a parallel filter, or a combination of the two to further enhance the quality of the preconditioner. The decision whether those additional steps are taken is based upon the required accuracy and can be freely selected, depending on user requirements. 

\section{Parallelization details and issues}
\label{section:paralleldetails}

The previous algorithm can be split into the following 5 phases (Notice that phases 1 and 5 are only necessary when the initial matrix is not a \emph{diagonally dominant matrix (ddm)}): 
1) Initial matrix is transformed into a ddm, 2) Transformation of ddm for suitable \emph{Neumann series expansion}, 3) Monte Carlo method is applied to calculate sparse approximation of the inverse matrix, 4) Given 2, calculate the inverse of the ddm  from 3, 5) Recovery process is applied to calculate the inverse of the original matrix due to the transformation in 1. 
It must be noted that the last phase requires in general $\mathcal{O}(n^3)$ operations and hence is generally neglected. Prior numerical experiments
have demonstrated that it is not compulsory to obtain an effective preconditioner.

This algorithm was originally designed for a HPC cluster composed of single-core compute nodes. It is written in C and uses the MPI library. It also makes use of the BeBOP sparse matrix converter~\cite{hoemmen2011bebop} to translate the input matrix format into a CSR format.

\subsection{MPI implementaion}\label{subsec:Parallelization-MPI}
Matrices $A$, $B_1$ and $P$ (the transition probability matrix) are calculated during the phases mentioned above. Note that $A = (I - C),C=B_1^{-1}\hat{B}$ and $B_1 = diag(\hat{B})$. 
Then a procedure is called by all the processes in which the partitioning of the matrix $A$ is carried out. The distribution of the work is done evenly when the number of rows is divisible by the number of processes. In the opposite case, the remaining rows are distributed among the smaller MPI processes (without including the Master process).
After that,  matrices $A$, $B_1$ and $P$ are \emph{broadcast} using MPI\_Broadcast(). Then the Monte Carlo process (phase 3) is started in parallel by all MPI processes.
During the Monte Carlo phase, each MPI process will calculate a piece of the inverse matrix of $C$ ($C^{-1}$), using matrix $A$; remember that $C = (I - A)$.
Column-scaling by $B_1^{-1}$ will then be applied to each row, to get the respective part of ${\hat{B}}^{-1}$ (phase 4).
After finishing the Monte Carlo process and phase 4, each process will send its part of the matrix ($\hat{B}^{-1}$) to the master process by calling MPI\_Send(). The master process will perform a corresponding MPI\_Receive() and will merge the received parts with its own.
Given a concatenation issue due to the CSR format, the Send-Receive process has to be ordered, having to receive first the data from process 1, then process 2 and so on.
Finally the last phase (5) is optionally executed by the master processes on matrix ($\hat{B}^{-1}$) to calculate $B^{-1}$. This step is optional and must be
enabled explicitly.
This process is difficult to parallelize due to its iterative nature. On the other hand, using an approach in which each iteration is executed in parallel, would imply a high increment in the communications given, that a synchronization would be required at each iteration.

\subsection{GPU implementaion}\label{subsec:Parallelization-GPU}
Regarding the GPU implementation it must be noted, that due to the irregular data access and comparatively short computation kernels
the method appears to be ill-suited for a GPU. Nevertheless a GPU can be used to accelerate the computation of the preconditioner using $(\text{MC})^2\text{MI}$
if care is taken to keep the GPU sufficiently busy.

If the requirement for a sparse inverse is abandoned the algorithm, with or without recovery, yields itself well to an implementation for GPUs.
This restricts the dimensionality of the matrices the algorithm is applicable to to those that fit entirely into the main memory of the accelerator device.
If a preconditioner is to be computed using $(\text{MC})^2\text{MI}$ on one or more GPUs for large, sparse, matrices a non-negligible amount
of overhead is introduced in order to ensure that only the most relevant entries of the inverse are retained for each row.
In a first implementation for NVIDIA GPUs the entire sparse inverse
was stored on the device, along with the necessary (preprocessed) matrix. Since the number of different entries visited by a chain
is not known a-priori the entire set of Markov Chains is simulated at once and used to fill a contiguous array corresponding to one
row of the approximate inverse. Afterwards only a prescribed number of entries largest in magnitude are retained for the sparse approximate inverse.
The rationale behind this is that if the inverse is itself considered as a Markov Chain only the
entries largest in magnitude will contribute significantly to its inverse (the original matrix).
As was the case with the previous implementation an extension to multiple GPUs is comparatively simple and has therefore been implemented.
\section{Algorithmic Modifications}\label{sec:AlgModifications}

The original code, provided by Diego D\'avila, was corrected to adhere to the MPI 3.0 standard and therefore be portable. This was crucial for performance analysis
on the testing system c.f. sec. \ref{sec:NumExp}. Furthermore a parallelizable pseudo-random number generator (PRNG) was used to replace the original generator, which was not suited for parallel environments.
\subsection{Matrix Reduction}
The computation of an approximate inverse using Markov Chain Monte Carlo (MCMC) requires
the knowledge of the whole state space - hence of the entire matrix $A$. The distribution of $A$
among the parallel workers becomes increasingly expensive with growing matrix size.

An obvious way to accelerate the method is to reduce the amount of data being transferred, i.e., to reduce the number of non-zero entries
of the matrix. Since the magnitude of the entries of $A$ signifies its importance in the MCMC simulation we decided to
drop a set percentage of the smallest entries of the matrix. This modifies the linear system and hence the correctness
of the approach had to be verified.

\subsection{Implementation Specifics - MPI}\label{subsec:ImplementationSpecifics_MPI}
As a first step the pseudo-random number generator used in the original version of the program
was replaced by TRNG\footnote{https://www.numbercrunch.de/trng}. This was necessary since the original code used the standard C
PRNG, which does not possess a sufficiently long period to guarantee statistic independence of the Markov chains for large matrices. Additionally it is not designed for parallel environments. Both flaws are rectified by using TRNG.

The amount of communications has already been reduced to an almost-minimum in the previous implementation of $(\text{MC})^2\text{MI}$.
In one iteration of the improvement of the code the broadcast of the transition probabilities (necessary for the MCMC simulation) was
eliminated. Instead these probabilities were computed by every worker from its knowledge of $A$.

Furthermore a minor non-conformity to the MPI standard was eliminated, which made the code reliant upon a specific
implementation of MPI, thus preventing the use of the preferred compiler and optimized MPI implementation on MareNostrum 4.
\subsection{Implementation Specifics - GPU}\label{subsec:ImplementationSpecifics_GPU}
Compared to the host machine the GPU has a very limited amount of memory and requires a more elaborate
approach to memory handling. Due to memory constraints storage of a dense block of an inverse on the device is not feasible,
and neither is on-the-fly transfer of computed entries to the host - due to latency constraints.
We have opted to allocate and fill a block of the sparse inverse on the device and transfer it to
the host matrix at the end of the computation.
This differs from the MPI implementation in so far as the computation of each row requires additional memory management overhead
but the final reduction of the separate blocks of the inverse is cheaper since the necessary storage and data layout is known beforehand. The downside being that
for some matrices entries of the inverse may be lost for some rows, whilst others contain unused entries ($=0$). This deficiency will be addressed in future 
versions of the GPU implementation.
A further difference from the MPI implementation is the usage of $\alpha\cdot\Vert B\Vert \cdot\text{sgn}(B_{i,i})$ as entries of the matrix $B_2$, as opposed
to $\alpha\cdot\Vert B\Vert$. This ensures that even if the signs of the diagonal elements are non-uniform the augmentation will yield a diagonally-dominant matrix.
This approach also reduces the perturbation of the original matrix caused by the augmentation procedure.
Usage of multiple GPUs was implemented by letting each device be controlled by a dedicated OpenMP process. 
\section{Numerical Experiments}\label{sec:NumExp}

	\begin{table}[h]
		 \caption{Matrix set.}
		\tabcolsep=0.1cm
		\begin{tabular}{c c c c c}
			\hline
			\textbf{Matrix} & \textbf{Dimension} & \textbf{Non-zeros} & \textbf{Sparsity} \\
			 \hline
			ID1\_2\_P3\_7\_stiffness & $514,369 \times 514,369$ & 8,702,911 & 0.003\% \\
			 nonsym\_r3\_a11 & $ 20,930 \times 20,930 $ & 638,733 & 0.15\% \\
			 rdb2048\_noL & $ 2,048 \times 2,048 $ & 12,032 & 0.29\% \\
			 sym\_r6\_a11 & $ 1,314,306 \times 1,314,306 $ & 36,951,316 & 0.02\% \\
			 \hline
		 \end{tabular}
     	 \label{tbl:matrix_set}   
	\end{table}
	
\subsection{Execution Environment}
The set of matrices chosen for the assessment of the proposed modifications is listed in tbl. \ref{tbl:matrix_set}.
The set contains symmetric and non-symmetric matrices of varying sizes and filling fractions. Matrices \verb+nonsym_r3_a11+
and \verb+sym_r6_a11+ have been provided by our collaborators and are representative
of systems occurring in climate simulations. The matrix \verb+rdb2048+ has been taken
from the Florida University's matrix collection and \verb+ID1_2_P3_7_stiffness+ is a discretized Laplacian using cubic finite elements on a fine mesh.

Almost all numerical experiments were carried out on the MareNostrum 4 (MN4) cluster at the Barcelona Supercomputing Centre in Spain.
The machine consists of 3456 nodes with 2 Intel Xeon Platinum 8160\makeatletter @\makeatother 2.1 GHz per node.
The nodes are connected via Intel Omni-Path HFI Silicon 100 (100 GBit/s) adapters.
The evaluation of the preconditioners was performed using 3 nodes of MareNostrum 4. The number was chosen arbitrary but kept constant,
thereby ensuring that the execution times of the preconditioned iterative solvers would be comparable for preconditioners computed
using CPUs and GPUs.

Earlier experiments evaluating the performance of Tesla K80 GPUs were performed on a GPU workstation and on 
the institutional cluster set up by AL at the Institute for Theoretical Physics.
Said cluster consists of 12 Nodes connected by a common 10GBit ethernet network and each containing two Intel Xeon E5-2640v4 CPUs.

On MN4 both, $(\text{MC})^2\text{MI}$ and MSPAI, were compiled using the INTEL compiler (v 17.0.4) and MPI implementation (build 20170405).
Execution was carried out in exclusive mode with CPU clock speeds fixed to the second-highest speed-step using batch script options to SLURM.

The computed preconditioners were validated using the GMRES implementation provided by Trilinos(v. 12.10.1).

For most experiments a precision of $\epsilon,\delta = 2^{-4}$ was chosen for MSPAI and MCMCMI. Additionally, for MCMCMI 
the scaling of the diagonal was performed using $\alpha = 5$.

To ensure that the GPUs are well-utilized a precision of $\epsilon,\delta \in\lbrace 0.01,0.005\rbrace$ has been chosen
in the numerical experiments comparing GPUs and CPUs. This choice provides a first limit on the range of parameters for which the use
of of $(\text{MC})^2\text{MI}$ on accelerators could be considered.
\subsection{Fitness of purpose}\label{subsec:FitnessForPurpose}

All of the numerical experiments in this section have been carried out with a fixed execution configuration of $48$ processes
spread evenly over two nodes of MN4.
\begin{figure}[h!]
\includegraphics[height=0.85\linewidth, width=0.95\linewidth]{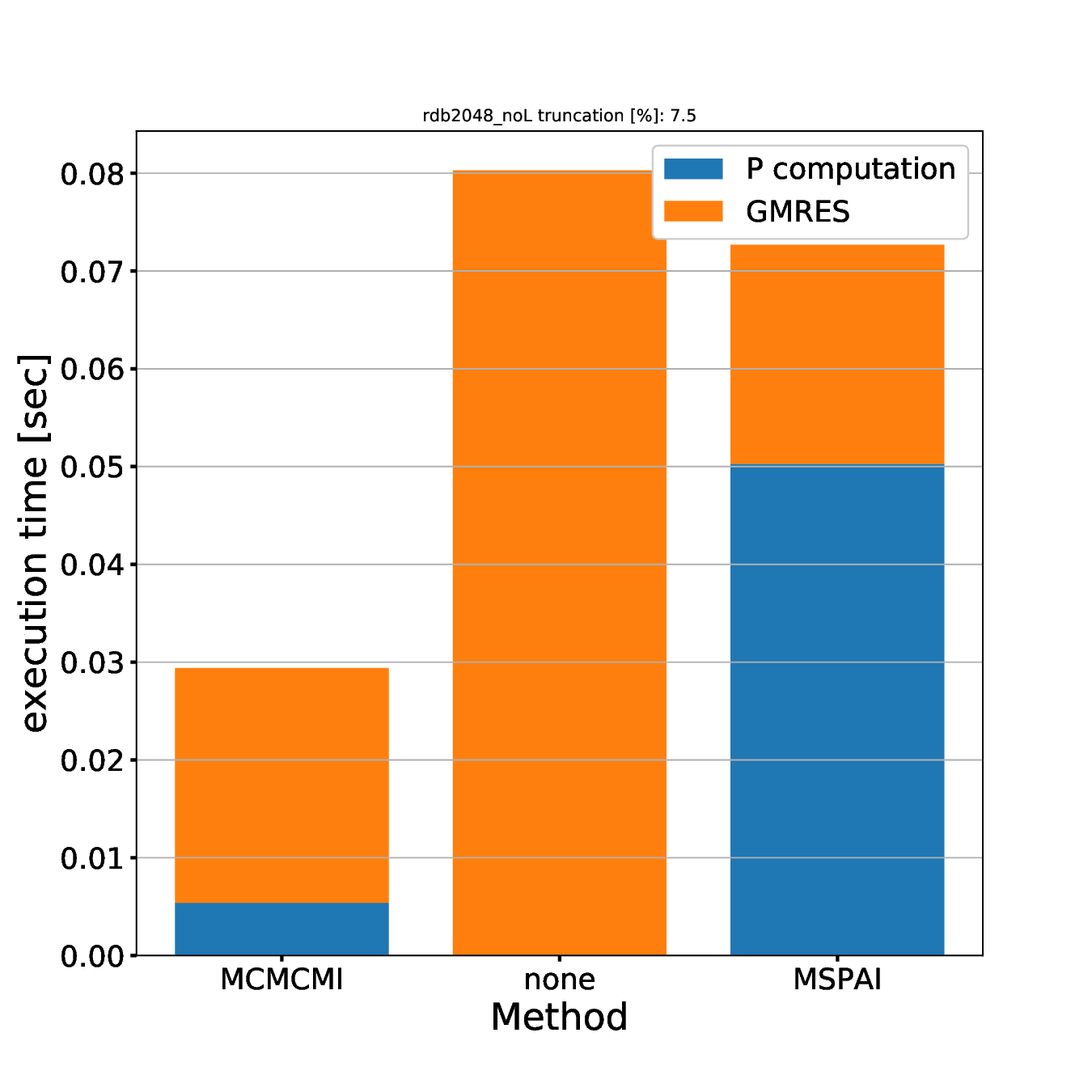}
\caption{Total execution time for \texttt{rdb2048\_noL} with $7.5\%$ of the value range of the entries removed.}
\label{fig:MC_vs_MSPAI_total_rdb}
\end{figure}

\begin{figure}[h!]
\includegraphics[height=0.85\linewidth, width=0.95\linewidth]{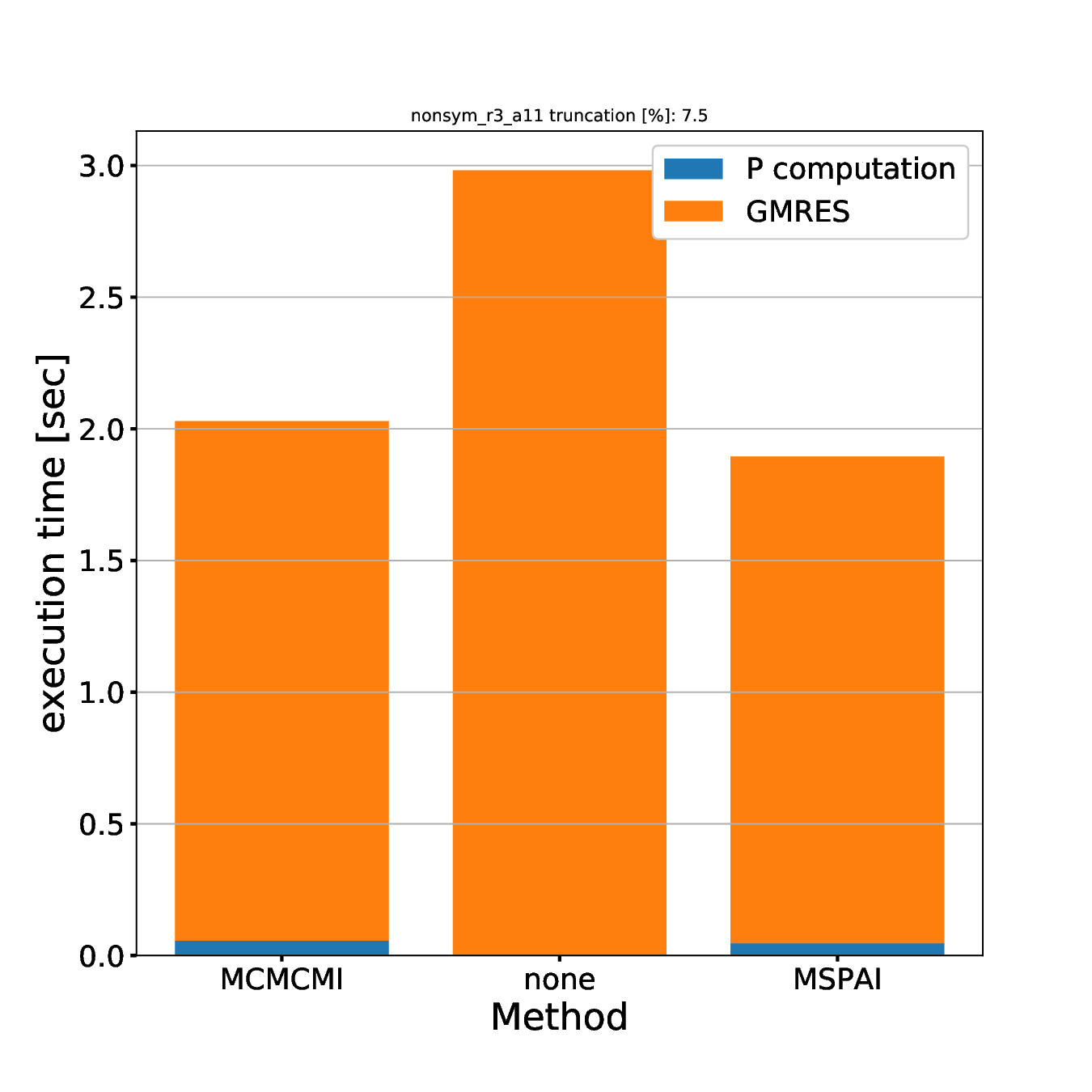}
\caption{Total execution time for \texttt{nonsym\_r3\_a11} with $7.5\%$ of the value range of the entries removed.}
\label{fig:MC_vs_MSPAI_total_nonsym}
\end{figure}

\begin{figure}[h!]
\includegraphics[height=0.85\linewidth, width=0.95\linewidth]{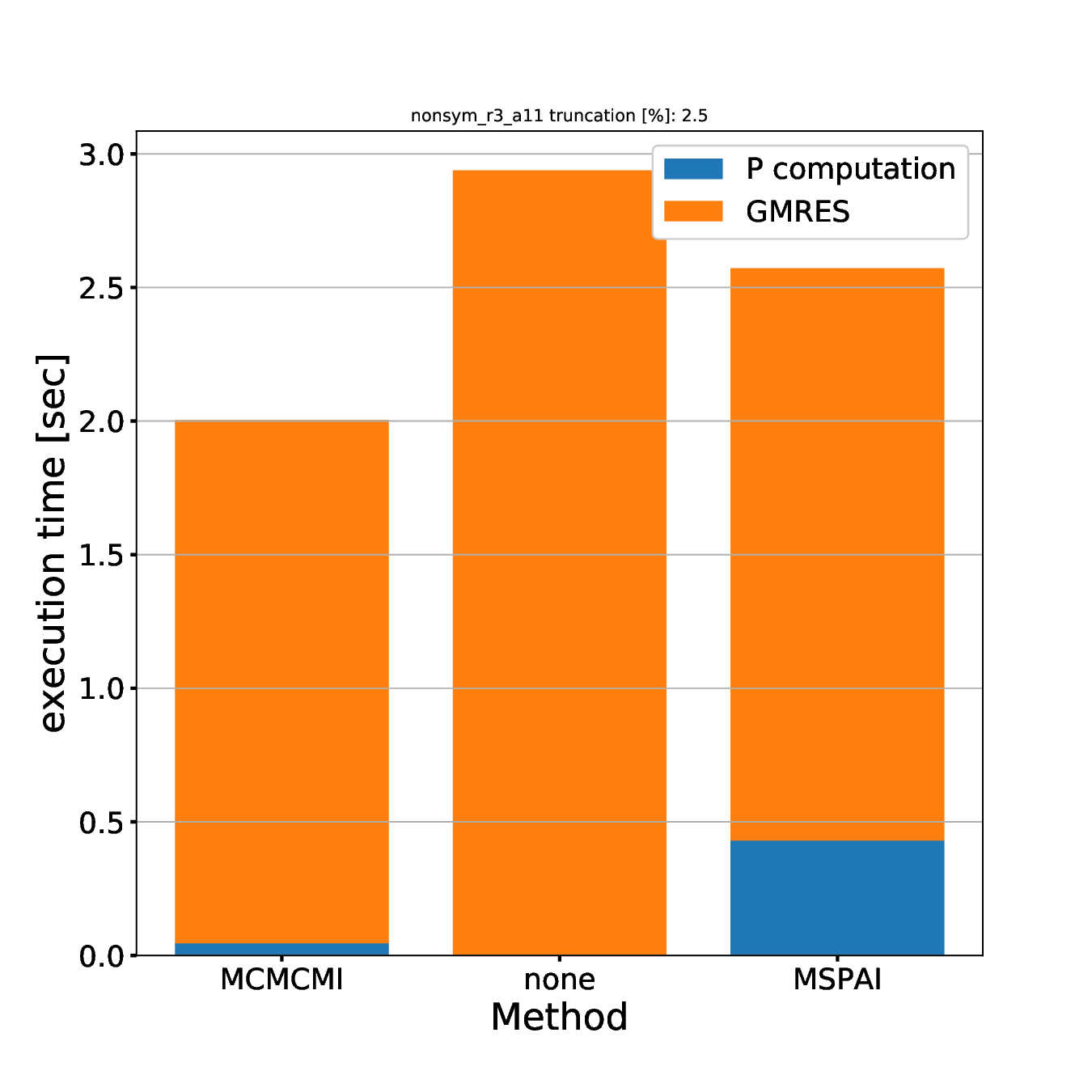}
\caption{Total execution time for \texttt{nonsym\_r3\_a11}. When using GMRES with a termination condition $\tfrac{\Vert r\Vert_2}{\Vert b\Vert_2}\leq 10^{-6}$
and a precision of $\epsilon = 0.0625$ for the computation of the preconditioner and removing $2.5\%$ of the entries smallest in magnitude.}
\label{fig:MC_vs_MSPAI_vs_Pure_RT_nonsymDrop}
\end{figure}

The total execution time (preconditioner computation and GMRES execution) is provided in fig. \ref{fig:MC_vs_MSPAI_total_rdb}
and \ref{fig:MC_vs_MSPAI_total_nonsym} for two different matrices.
Henceforth \verb+none+ refers to the method without preconditioner and the preconditioner (computed using MSPAI or MCMCMI) is designated $P$.
MSPAI is more effective in the case of the larger of the two matrices but
only in the case when $>7 \%$ of the value range of the elements of the matrix have been dropped. If less elements are removed (c.f. fig.
 \ref{fig:MC_vs_MSPAI_vs_Pure_RT_nonsymDrop}) $(\text{MC})^2\text{MI}$ will require less computation time.

\begin{figure}[h!]
\includegraphics[height=0.85\linewidth, width=0.95\linewidth]{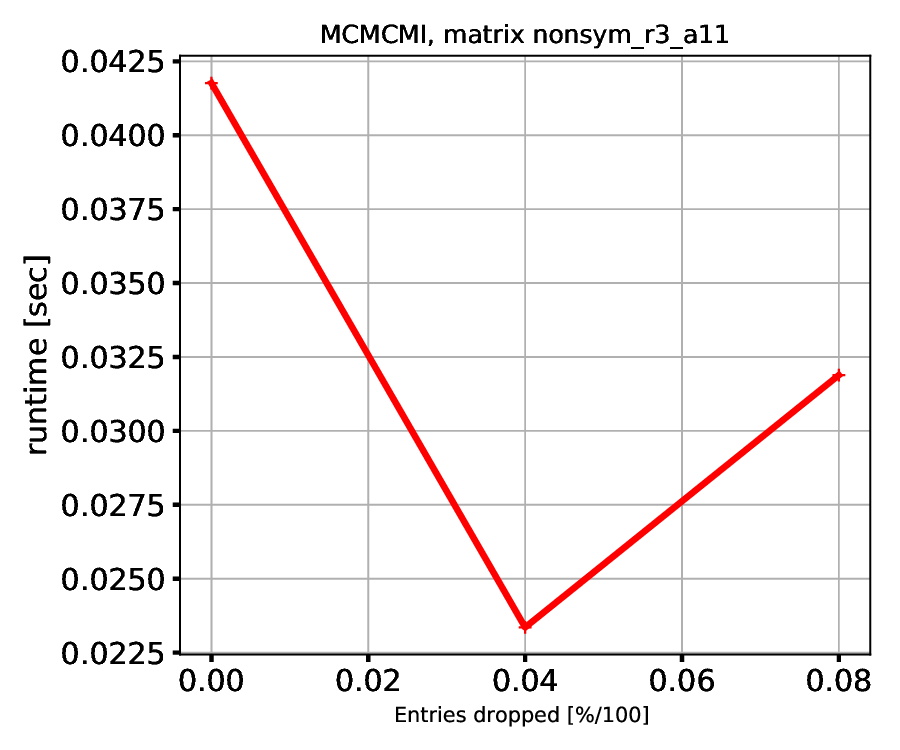}
\caption{Execution time of the preconditioner computation.}
\label{fig:MC_RT_vs_Drop}
\end{figure}

In fig. \ref{fig:MC_RT_vs_Drop} one can see, that the idea of removing a set amount of small elements may well accelerate the
computation of the preconditioner. The outcome depends on the matrix and there will, in general, be an optimal amount of negligible
entries for each matrix. In the case of \verb+nonsym_r3_a11+ that amount is between $2\%$ and $6\%$. If more entries 
are removed the amount of information contained in the matrix becomes insufficient to create a good preconditioner.

As can be seen in fig. \ref{fig:MC_vs_MSPAI_vs_Pure_RT_nonsymDrop} the reduction of the amount of information required to be broadcast, coupled with a
moderate precision requirements for the approximate inverse will result in a shorter overall execution time when a preconditioner is computed using $(\text{MC})^2\text{MI}$.
This demonstrates that the method may be used in cases where the preconditioner has to be recomputed every time prior to its usage (i.e.,
in iterative methods where the matrix changes in every step).

Finally we attempted to use the Monte Carlo method to compute a preconditioner for the \texttt{bcsstk38} matrix of the sparse matrix collection,
whose condition number surpasses $5\cdot 10^{16}$ and which has a non-trivial nullity.
Accordingly the iterative solver used to test the preconditioner for this case (BiCGstab) fails to converge
if no preconditioner is used, reaching the defined upper bound of $30000$ iterations for a desired precision of $\tfrac{\Vert r\Vert_2}{\Vert b\Vert_2}\leq 0.45$.
Using the preconditioner computed with $(\text{MC})^2\text{MI}$ for $\epsilon=0.01$ enables BiCGstab to converge,
reducing the number of steps required to achieve the desired bound to $3852$
and the total execution time (preconditioner + BiCGstab) from $9.7[sec]$ to $1.3[sec]$ - in this case using $96$ instead of $48$ processes.

\subsection{Scaling to Moderate Number of Cores/Processors}
\begin{figure}[htb]
\includegraphics[height=0.85\linewidth, width=0.95\linewidth]{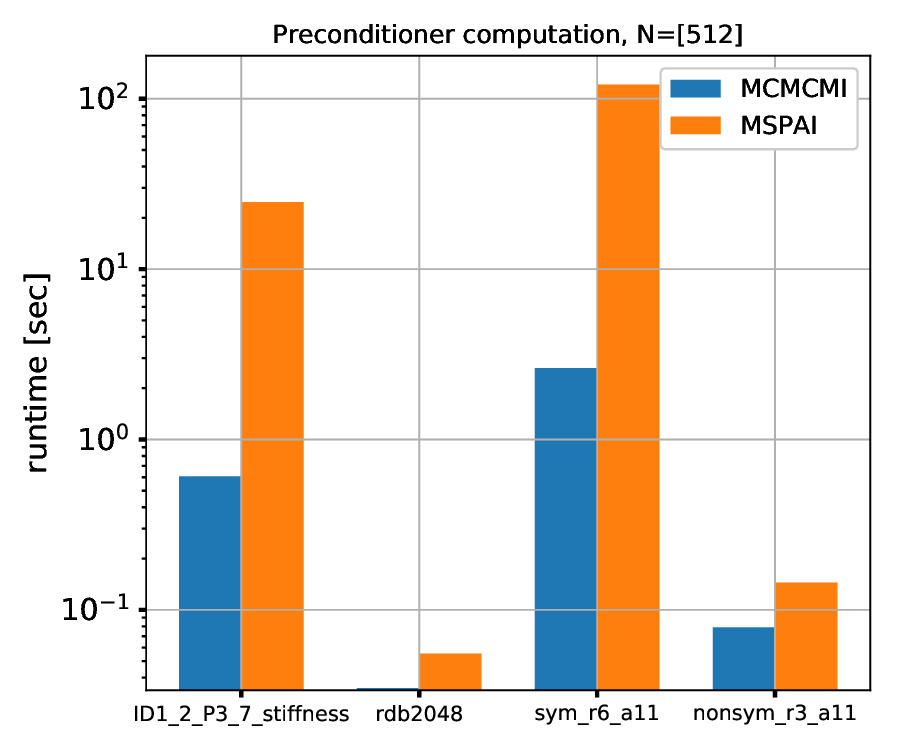}
\caption{Execution time of the preconditioner computation for different matrices. Transition probabilities were computed 
by the master process and broadcast to the workers.}
\label{fig:MC_vs_MSPAI_runtime_vs_matrix}
\end{figure}
In fig. \ref{fig:MC_vs_MSPAI_runtime_vs_matrix} the execution time of the preconditioner computation using $512$ processes
is shown for all matrices of tbl. \ref{tbl:matrix_set}. It is obvious that $(\text{MC})^2\text{MI}$ is
superior to MSPAI in every case, with the largest savings being achieved for rather dense or very large matrices.
Note that this is purely a comparison of the time required to compute a preconditioner using the appropriate method.
\subsection{$(\text{MC})^2\text{MI}$ on Accelerators}
\begin{figure}[htb]
\includegraphics[height=0.85\linewidth, width=0.95\linewidth]{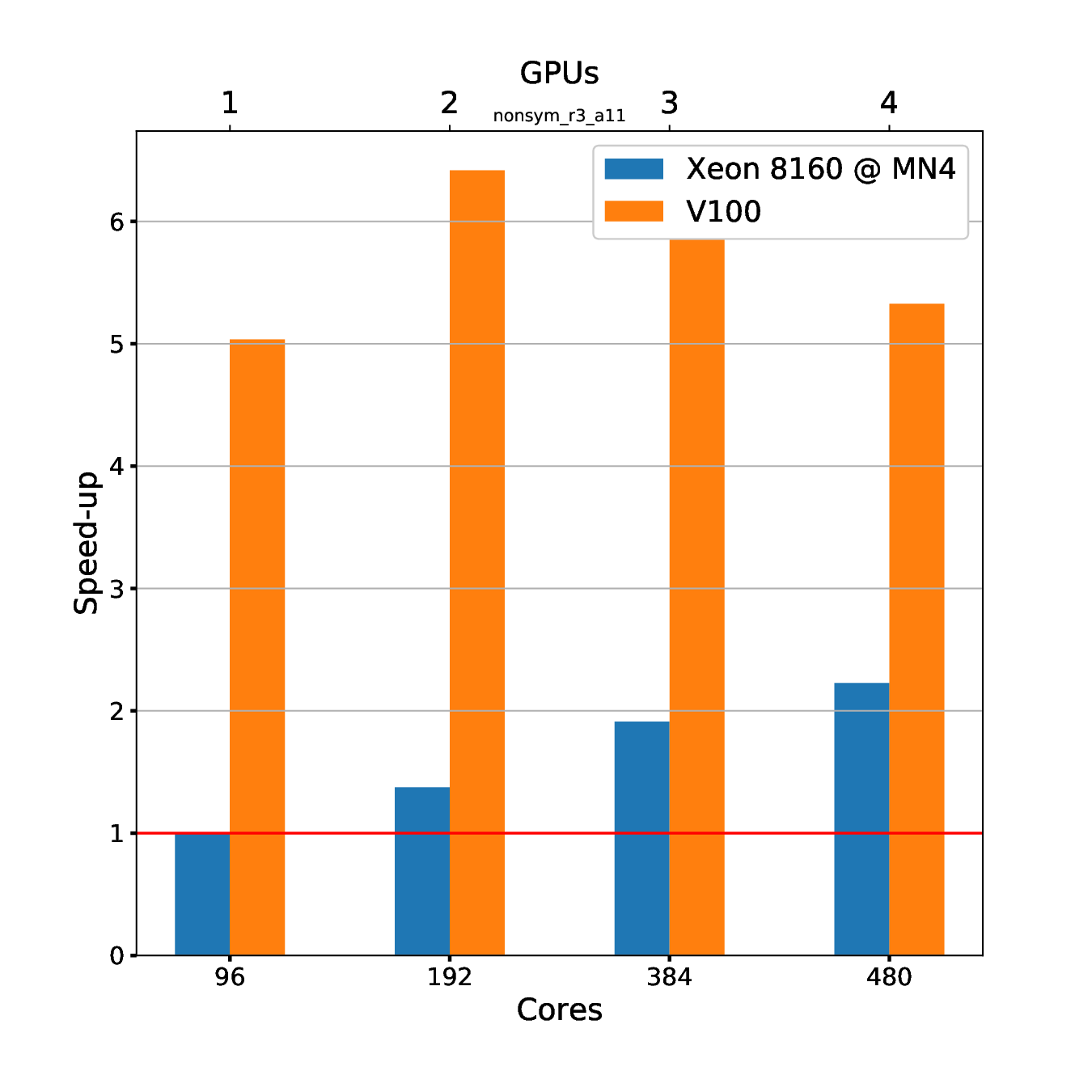}
\caption{Speed-up of the calculation performed on a NVIDIA V100 in comparison to CPU cores for the nonsym\_r3\_a11 matrix
and a precision of $\epsilon = 0.01$. The red line represents $1\times$.}
\label{fig:nonsym_GPU_speedup}
\end{figure}
In a final set of numerical experiments we investigated the feasibility of using accelerators, specifically NVIDIA GPUs 
to speed up the computation of the preconditioners using $(\text{MC})^2\text{MI}$. To this end the algorithm delineated in sec.
\ref{subsec:Algorithm} was implemented in CUDA and evaluated for matrices of tbl. \ref{tbl:matrix_set}. 
Here we have to note that unlike for the pure MPI implementation a sufficiently small $\epsilon,\delta$ (i.e., a high precision)
is necessary to fully utilise the GPU as such a precision of $\epsilon = 0.01$ was chosen for the experiments.
The latter were performed on Tesla K80  as well as Volta V100 devices using a variable number of GPUs.

Fig. \ref{fig:nonsym_GPU_speedup} shows the typical behaviour of the Markov Chain Monte-Carlo method when implemented on GPUs
using \texttt{nonsym\_r3\_a11} as an example. The execution time of the pure MPI implementation on two nodes of MareNostrum 4 serves as 
a reference. It is immediately obvious that a GPU is significantly faster by up to a factor of $\sim 6.5$.
The speed-up decreases when using 3 or more GPUs, which is to be attributed to the overhead introduced by the memory management.
Profiling results indicate that in this case the time required to sort the entries of the inverse row matches the time required to compute them
using $(\text{MC})^2\text{MI}$. An additional factor limiting the performance, which has not yet been eliminated, is the necessity to compact
the pre-processed matrix on the host before the MC iteration may be performed. In the present case this reduces the achievable speed-up by a factor of $\sim 2$.
Note that for this comparison the $\alpha$ parameter was chosen to be $4.0$ instead of $5.0$. This change results in a $~20\%$ longer execution time.
\begin{figure}[htb]
\includegraphics[height=0.85\linewidth, width=0.95\linewidth]{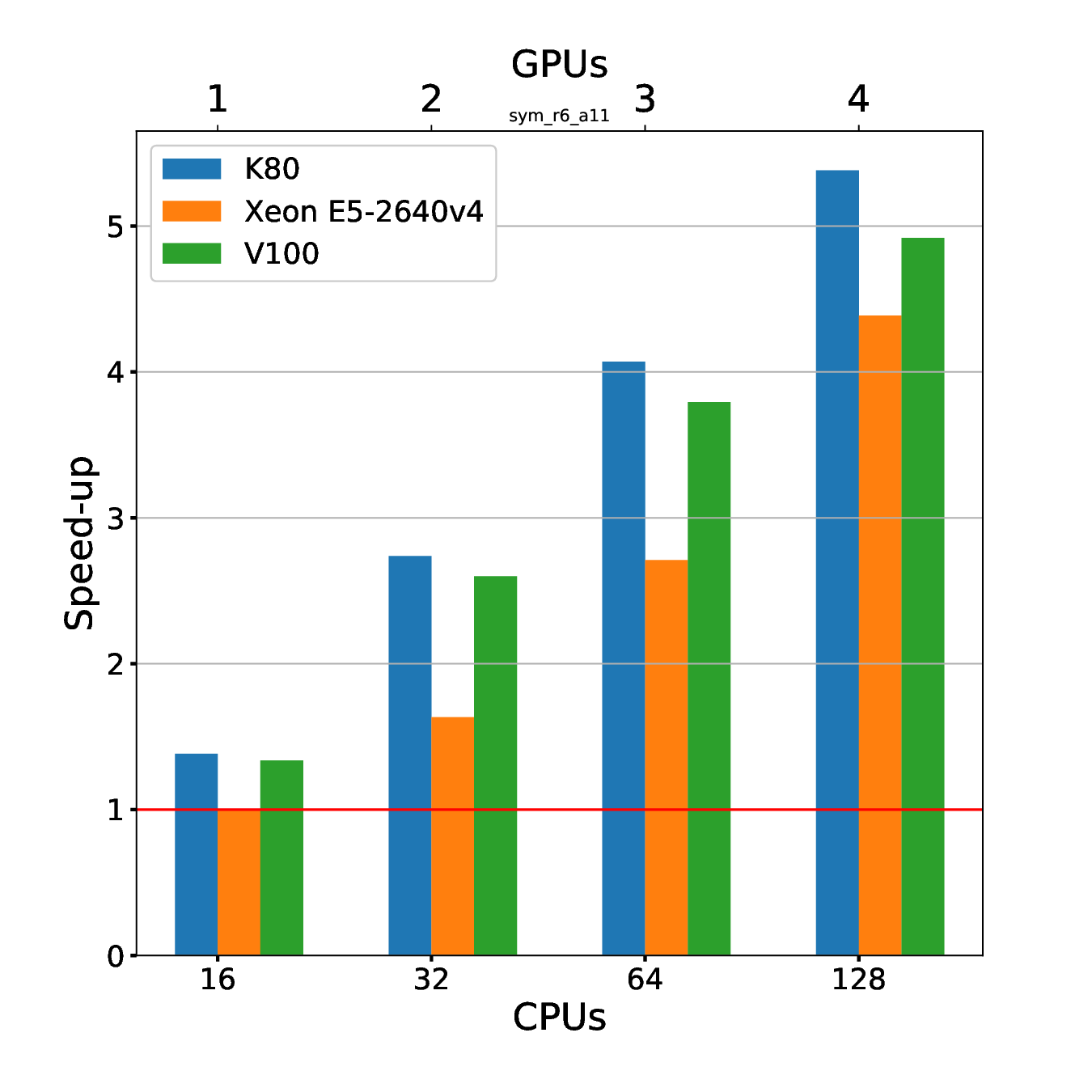}
\caption{Speed-up of the calculation performed on a NVIDIA V100 in comparison to CPU cores of a Broadwell cluster and to NVIDIA
K80 GPUs for the sym\_r6\_a11 matrix. The red line represents $1\times $.}
\label{fig:sym_V100_vs_K80_vs_GPU_speedup}
\end{figure}
The effect of an increased amount of work can be seen in fig. \ref{fig:sym_V100_vs_K80_vs_GPU_speedup}, where the speed-up achieved
in comparison to an older CPU architecture is shown. The comparison is provided due to the given CPU and GPU resources being 
an easily accessible resource maintained by AL at the Institute for Theoretical Physics in T\"{u}bingen as well as to their availability
to common users (in comparison to multiple V100).
The speed-up provided by the older Teslas is limited for the given case due to the sparsity of the matrix. Further numerical
experiments indicate that utilisation of the GPU may be improved by increasing the desired precision.

\begin{figure}[htb]
\includegraphics[height=0.85\linewidth, width=0.95\linewidth]{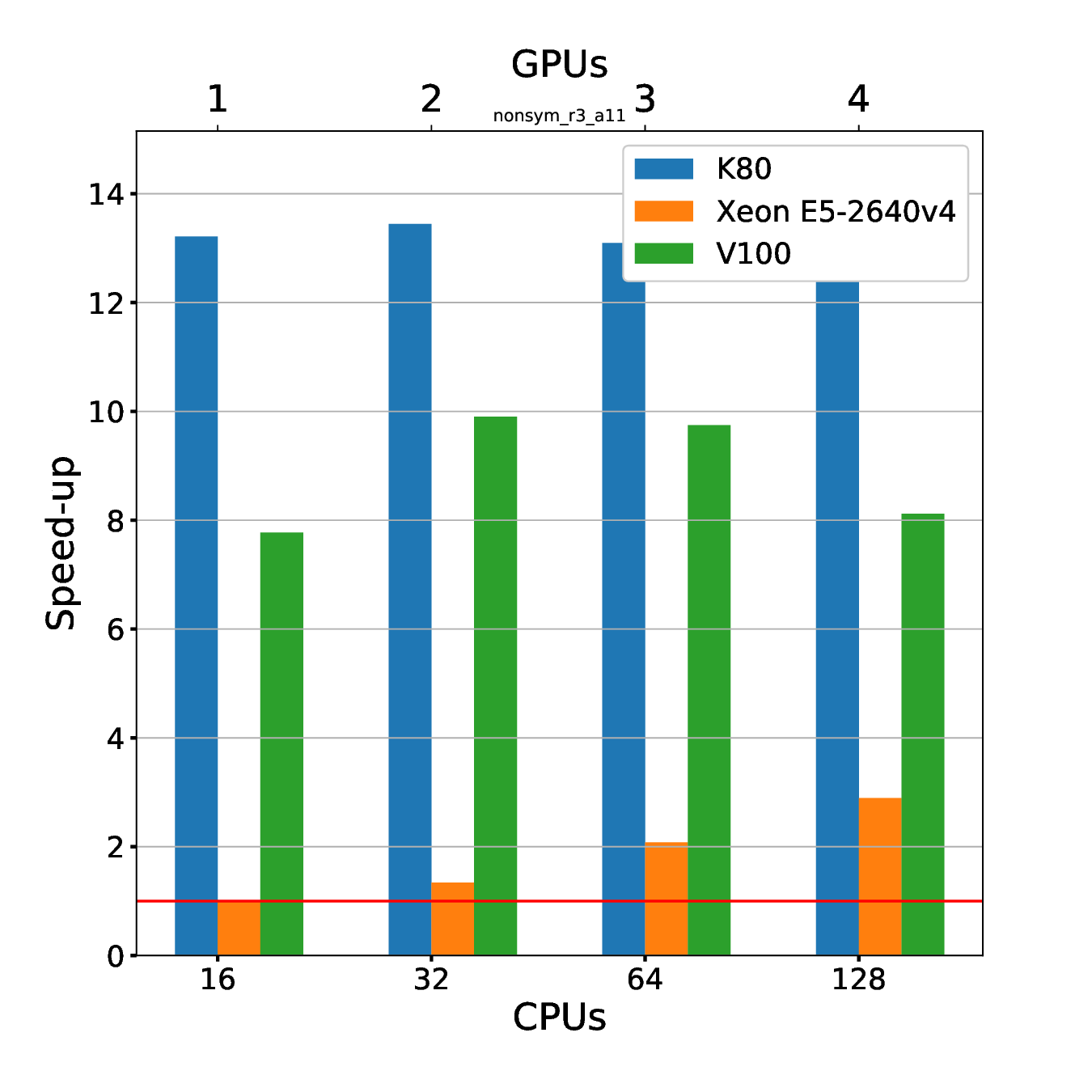}
\caption{Speed-up of the calculation performed on a NVIDIA V100 in comparison to CPU cores of a Broadwell cluster and to NVIDIA
K80 GPUs for the nonsym\_r3\_a11 matrix. The red line represents $1\times $.}
\label{fig:nonsym_V100_vs_K80_vs_GPU_speedup}
\end{figure}
Finally fig. \ref{fig:nonsym_V100_vs_K80_vs_GPU_speedup} shows the speed-up achieved by two generations of NVIDIA GPUs for the nonsym\_r3\_a11 matrix
compared to the small institutional cluster in T\"{u}bingen. As has been demonstrated in fig. \ref{fig:nonsym_GPU_speedup} the amount of work provided by this matrix is insufficient to mask the overhead of data management and the CPU portion of the preprocessing stage. Both are currently being adressed in development.
The striking feature is that the newer architecture appears to perform worse than the older one. We believe this to be an artefact due to
an insufficient optimization of the GPU code for the NVIDIA Volta architecture, since it has been originally developed and optimized for the Kepler
architecture.

\section{Conclusions and Future Work}
\label{section_conclusion}
In summary we have shown that the computation of a preconditioner using $(\text{MC})^2\text{MI}$ method can be accelerated
by ballancing the precision with which the preconditioner is calculated as well as dropping entries of the original matrix depending on the precision. The quality of the resulting preconditioner
does not deteriorate as fast as is the case if the same approach is applied to MSPAI. 
The approach shows that in most cases the number of iterations required by GMRES or BiCGstab to solve the resulting system of Linear Algebraic Equations can be substantially reduced. If only a rough estimate of the inverse is required the combination of $(\text{MC})^2\text{MI}$ and an appropriate (for the matrix type)
iterative method can result in a lower total execution time, when compared to a non-preconditioned method.

The numerical experiments indicate that for $\epsilon,\delta < 0.01$ (at high precisions) the usage of GPUs should be considered.
It has been demonstrated, that despite the apparent bad suitability for a GPU the $(\text{MC})^2\text{MI}$ method may still be successfully
used with it.
Future work will focus on a merging of the CPU and GPU implementations using the tasking constructs of OpenMP 4.5. This approach
promises to reduce the overhead of memory management on the GPU whilst simultaneously utilising the host to its full extend.
Preliminary profiling suggests a potential increase in performance by a factor of $\gtrsim 2$.
Furthermore an integrated application test for the Markov Chain Monte Carlo preconditioners is planned, to observe the performance 
on a wider set of matrices than the set used so far, as well as an investigation of the potential for a MPI+CUDA parallelisation of the method.
On the host side the pure MPI implementation will be rewritten to utilise hybrid parallelism using MPI+OpenMP and implement a better load balancing.

\section{Acknowledgments}
Anton Lebedev wishes to thank the Severo Ochoa program, Spain,  for providing a mobility grant enabling him to work on this project at the Barcelona Supercomputing Centre.

\bibliographystyle{IEEEtran}
\bibliography{references}

\begin{thebibliography}{10}
\providecommand{\url}[1]{#1}
\csname url@samestyle\endcsname
\providecommand{\newblock}{\relax}
\providecommand{\bibinfo}[2]{#2}
\providecommand{\BIBentrySTDinterwordspacing}{\spaceskip=0pt\relax}
\providecommand{\BIBentryALTinterwordstretchfactor}{4}
\providecommand{\BIBentryALTinterwordspacing}{\spaceskip=\fontdimen2\font plus
\BIBentryALTinterwordstretchfactor\fontdimen3\font minus
  \fontdimen4\font\relax}
\providecommand{\BIBforeignlanguage}[2]{{%
\expandafter\ifx\csname l@#1\endcsname\relax
\typeout{** WARNING: IEEEtran.bst: No hyphenation pattern has been}%
\typeout{** loaded for the language `#1'. Using the pattern for}%
\typeout{** the default language instead.}%
\else
\language=\csname l@#1\endcsname
\fi
#2}}
\providecommand{\BIBdecl}{\relax}
\BIBdecl

\bibitem{golub1996matrix}
\BIBentryALTinterwordspacing
G.~Golub and C.~Loan, \emph{Matrix computations}, ser. Johns Hopkins studies in
  the mathematical sciences.\hskip 1em plus 0.5em minus 0.4em\relax Johns
  Hopkins University Press, 1996. [Online]. Available:
  \url{http://books.google.es/books?id=mlOa7wPX6OYC}
\BIBentrySTDinterwordspacing

\bibitem{AlexandrovStrassburgA14}
\BIBentryALTinterwordspacing
J.~Stra{\ss}burg and V.~N. Alexandrov, ``Enhancing monte carlo preconditioning
  methods for matrix computations,'' in \emph{Proceedings of the International
  Conference on Computational Science, {ICCS} 2014, Cairns, Queensland,
  Australia, 10-12 June, 2014}, 2014, pp. 1580--1589. [Online]. Available:
  \url{http://dx.doi.org/10.1016/j.procs.2014.05.143}
\BIBentrySTDinterwordspacing

\bibitem{AlexandrovV15}
\BIBentryALTinterwordspacing
V.~N. Alexandrov and O.~A. Esquivel{-}Flores, ``Towards monte carlo
  preconditioning approach and hybrid monte carlo algorithms for matrix
  computations,'' \emph{Computers {\&} Mathematics with Applications}, vol.~70,
  no.~11, pp. 2709--2718, 2015. [Online]. Available:
  \url{http://dx.doi.org/10.1016/j.camwa.2015.08.035}
\BIBentrySTDinterwordspacing

\bibitem{vajargah2001parallel}
F.~Vajargah \emph{et~al.}, ``Parallel monte carlo algorithms for matrix
  computations,'' Ph.D. dissertation, University of Reading, 2001.

\bibitem{strassburgPhd2014}
J.~Strassburg, ``{On hybrid and resilient Monte Carlo methods for linear
  algebra problems},'' Ph.D. dissertation, The University of Reading, 2014.

\bibitem{branford2008}
S.~Branford, ``{Hybrid Monte Carlo Methods for Linear Algebra Problems},''
  Ph.D. dissertation, School of Systems Engineering, The University of Reading,
  April 2009.

\bibitem{carpentieri2000some}
B.~Carpentieri, I.~Duff, and L.~Giraud, ``{Some sparse pattern selection
  strategies for robust Frobenius norm minimization preconditioners in
  electromagnetism},'' \emph{Numer. Linear Algebra Appl}, vol.~7, pp. 667--685,
  2000.

\bibitem{carpentieri2000experiments}
B.~Carpentieri, L.~Giraud \emph{et~al.}, ``Experiments with sparse
  preconditioning of dense problems from electromagnetic applications,''
  CERFACS, Toulouse, France, Tech. Rep., 2000.

\bibitem{alleon1997sparse}
G.~All{\'e}on, M.~Benzi, and L.~Giraud, ``Sparse approximate inverse
  preconditioning for dense linear systems arising in computational
  electromagnetics,'' \emph{Numerical Algorithms}, vol.~16, no.~1, pp. 1--15,
  1997.

\bibitem{Branford2008grid}
\BIBentryALTinterwordspacing
S.~Branford, C.~Sahin, A.~Thandavan, C.~Weihrauch, V.~Alexandrov, and I.~Dimov,
  ``{Monte Carlo Methods for Matrix Computations on the Grid},'' \emph{Future
  Generation Computer Systems}, vol.~24, no.~6, pp. 605 -- 612, 2008. [Online].
  Available:
  \url{http://www.sciencedirect.com/science/article/B6V06-4P8GWNW-1/2/a4fbfc616f7b37f8237a0e484563188f}
\BIBentrySTDinterwordspacing

\bibitem{Dimov1998convergent}
I.~Dimov and V.~Alexandrov, ``{A New Highly Convergent Monte Carlo Method for
  Matrix Computations},'' \emph{Mathematics and Computers in Simulation},
  vol.~47, no. 2-5, pp. 165--181, Aug 1998.

\bibitem{Dimov2007}
I.~Dimov, V.~Alexandrov, R.~Papancheva, and C.~Weihrauch, ``{Monte Carlo
  Numerical Treatment of Large Linear Algebra Problems},'' in \emph{Lecture
  Notes in Computing Sciences: Computational Science - ICCS 2007}, vol.
  4487.\hskip 1em plus 0.5em minus 0.4em\relax Berlin: Springer-Verlag GmbH,
  2007, pp. 747--754.

\bibitem{mcrex2013}
T.~Evans, S.~Hamilton, W.~Joubert, and C.~Engelmann, ``{MCREX - Monte Carlo
  Resilient Exascale Project},'' http://www.csm.ornl.gov/newsite/documents/ \\
  CSMDSummer2013Newsletter.pdf, 2013.

\bibitem{benzi1996sparse}
M.~Benzi, C.~Meyer, and M.~Tuma, ``{A Sparse Approximate Inverse preconditioner
  for the Conjugate Gradient Method},'' \emph{SIAM Journal on Scientific
  Computing}, vol.~17, no.~5, pp. 1135--1149, 1996.

\bibitem{Huckle2010}
\BIBentryALTinterwordspacing
T.~Huckle, A.~Kallischko, A.~Roy, M.~Sedlacek, and T.~Weinzierl, ``{An
  efficient parallel implementation of the MSPAI preconditioner},''
  \emph{Parallel Computing}, vol.~36, no. 5-6, pp. 273 -- 284, 2010, parallel
  Matrix Algorithms and Applications. [Online]. Available:
  \url{http://www.sciencedirect.com/science/article/pii/S016781910900129X}
\BIBentrySTDinterwordspacing

\bibitem{grote2006spai}
M.~Grote and M.~Hagemann, ``{SPAI: SParse Approximate Inverse
  Preconditioner},'' \emph{Spaidoc. pdf paper in the SPAI}, vol.~3, p.~1, 2006.

\bibitem{huckle2003factorized}
T.~Huckle, ``Factorized sparse approximate inverses for preconditioning,''
  \emph{The Journal of Supercomputing}, vol.~25, no.~2, pp. 109--117, 2003.

\bibitem{strassburg2013scalability}
J.~Strassburg and V.~Alexandrov, ``{On Scalability Behaviour of Monte Carlo
  Sparse Approximate Inverse for Matrix Computations},'' in \emph{Proceedings
  of the ScalA 2013 Workshop}.\hskip 1em plus 0.5em minus 0.4em\relax ACM,
  2013, p.~6.

\bibitem{Fathi07}
\BIBentryALTinterwordspacing
B.~F. Vajargah, ``A new algorithm with maximal rate convergence to obtain
  inverse matrix,'' \emph{Applied Mathematics and Computation}, vol. 191,
  no.~1, pp. 280--286, 2007. [Online]. Available:
  \url{http://dx.doi.org/10.1016/j.amc.2007.02.085}
\BIBentrySTDinterwordspacing

\bibitem{hoemmen2011bebop}
M.~Hoemmen, R.~Vuduc, and R.~Nishtala, ``Bebop sparse matrix converter,''
  \emph{University of California at Berkeley. Web}, 2011.

\end{thebibliography}

\pagebreak

\appendices

\section{Artifact Description Appendix: On Advanced Monte Carlo Methods for Linear Algebra on Advanced Accelerator Architectures}

\subsection{Abstract}

We present observations on the performance of the implementation of the Markov Chain Matrix Inversion method for different versions of the x86
CPU architecture (Broadwell,Skylake) and NVIDIA GPUs of the Kepler and Volta architectures. The performance and correctness of the method
as a means of obtaining preconditioners for iterative systems is evaluated using Trilinos and compared to MSPAI.

The CPU (MPI) and GPU implementations of the Markov Chain method are compared to each other to determine the feasibility and limitations of
a GPU implementation of the method.

\subsection{Description}

\subsubsection{Check-list (artifact meta information)}

{\small
\begin{itemize}
  \item {\bf Algorithm: } Markov Chain Monte Carlo Matrix Inversion, preallocated row storage on CPU, stream compaction and sorting on GPU
  \item {\bf Compilation: }\emph{MareNostrum 4}: INTEL toolchain v 2017.4, with \verb+-xHost -O3+ optimization flags.
   \emph{ITP T\"{u}bingen: } GCC v. 6.3.1 (20170216)  with \verb+-march=native -O3 -mfma+ \verb+-malign-data=cacheline+ optimization flags.
  \item {\bf Run-time environment: }\emph{MareNostrum 4: }SLES 12-SP2, Kernel: 4.4.120-92.70.
  \emph{ITP T\"{u}bingen: }CentOS 7, Kernel: 3.10.0-514.26.2 (no KPTI mitigation).
  \emph{CTE Power} RHEL 7.4, Kernel: 4.11.0-44
  \item {\bf Hardware: }
  \emph{MareNostrum 4} Nodes with 2 Xeon Platinum 8160 CPUs each. 96GB RAM per node. Connected via 100GBit Intel Omni-Path HFI Silicon 100 in a fat tree network topology.
  \emph{CTE Power }(V100 machine) Nodes with 2 x IBM Power9 8335-GTG \@ 3.00GHz each. 512GB RAM and 4 V100 GPUs with 16GB HBM2 VRAM.
  \emph{ITP T\"{u}bingen }Nodes with 2 Intel Xeon E5-2640v4 CPUs each. 128GB RAM per node. Connected via 10GBit Ethernet, star network topology. Network parameters not optimized.
  \item {\bf Execution: } Via SLURM scheduler.
  \item {\bf Output: } Execution times (in milliseconds) are printed to standard output and processed from there.
  \item {\bf Experiment workflow: }Automated filling of SLURM script templates and automated enqueueing of the jobs by a generator script written in Python.
An index of numerical experiments is stored in the top-level directory where the generator script was called. This index is used to
collect and pre-process the results using an evaluation script written in Python.
Graphical analysis of the data is performed using a Jupyter notebook.
  \item {\bf Experiment customization: }Execution configuration of the job scripts customized to stay within storage quota.
  K80 experiments driven by a separate script.
  \item {\bf Publicly available?: }Currently not publicly available. Access to the authors repository can be granted upon request.
\end{itemize}
}

\subsubsection{How software can be obtained}
The GPU implementation can be obtained through the authors private Bitbucket repository upon request. The CPU implementation will be 
publicly available from said repository by the end of November.

\subsubsection{Hardware dependencies}
The optimal block and grid size of the GPU implementation are dependent on the used GPU and have hence to be adapted accordingly.
A rough search for minimal execution time using the sym\_r6\_a11 matrix suggested a block size of $96$ threads (3 warps) and a grid size of $170$ for the Volta GPUs.

\subsubsection{Software dependencies}

\paragraph*{CPU}
The CPU implementation uses version 4.15 of Tina's Random Number Generator library. The library implements parallel pseudorandom number generators
and is therefore key to the correctness of the presented method. It is available from www.numbercrunch.de/trng.
It also relies on the BeBOP sparse matrix library to handle CRS matrices.

\paragraph*{GPU}
The GPU verision has been implemented in C++ and CUDA. Of the CUDA libraries it utilises cuRAND in the core routines and
cuBLAS in some auxiliary routines and for testing purposes.
The V100 compilation was performed using CUDA Toolkit v9.1, the K80 compilation was performed using CUDA Toolkit v8.0.
Parsing of execution parameters is done using BOOST program options library (tested with BOOST $1.\lbrace 56,64,66\rbrace$) and
the Eigen linear algebra template library (http://http://eigen.tuxfamily.org/) to handle sparse matrices with a minor 
correction in the unsupported \verb+saveMarket+ routine.

\paragraph*{Testing}
Correctness checks of the preconditioners are carried out using a parallel implementation of CG/CGS/BiCG(stab)/GMRES.
The code performing these checks has been written in C++ and uses Trilinos (v 12.10.1 on MN4, 12.13 on the ITP cluster).

On MareNostrum 4 both the CPU implementation and the preconditioner testing code rely on the MPI implementation provided by INTEL.
On the ITP cluster the MPI library is MPICH 3.2.1.
\subsection{Installation}

\subsubsection{MareNostrum 4}
The MPI implementation is compiled using a simple Makefile and utilising the INTEL compiler \verb+mpiicc+
to compile all but the TRNG files, which are compiled using \verb+mpiicpc+. Compiler options are
\verb|-O3 -xHost -DPRECISION=1 -DP_UMMAO=0|
\verb|-DVER=3.0|
and linked to BeBOP libraries via
\verb|-Wl,-rpath=$(USRLIB) -lbebop_util|
\verb|-lsparse_matrix_converter|
and statically linked to the TRNG library \verb|libtrng4.a|
\subsection{CTE-POWER}
The GPU code is compiled using \verb+nvcc+ v9.1.85 with the following compiler flags\\
\verb|-std=c++11 -m64|
\verb| -Wno-deprecated-declarations|
\verb|-D EIGEN_NO_CUDA -arch=compute_70 -rdc=true|
\verb|-DNDEBUG -O3|
using a simple makefile which constitutes just a collection of source files to be compiled and linked.

\subsubsection{ITP T\"{u}bingen}
The process is the same as for the other two, except for the optimization flags:
\verb|-O3 -march=broadwell -ftree-vectorize|
\verb|-funroll-loops -ffunction-sections|
\verb|-malign-data=cacheline|

\subsection{Experiment workflow}

The numerical experiments carried out on the MareNostrum 4 and CTE-Power clusters at the Barcelona Supercomputing Centre were
executed in two stages:

\begin{enumerate}
\item \textbf{Stage}: Generate a set of preconditioners\\
The numerical experiments were executed using the SLURM scheduler. A generator script was written in Python. Said script accepts
a set of template files for the preconditioner computation and testing parameters as well as job scripts for generation and testing of preconditioners.
The execution parameters are collected in a separate parameter file and indexed by matrix in dictionaries.
The user may provide a desired number of repetitions the experiments will be run (10 as a default) -each repetition will
generate a preconditioner which will be stored with a file name containing the repetition number.
The generator script generates a directory structure and an index file for the desired numerical experiments. All of the jobs 
to generate preconditioners are launched using a simple launcher script and the generated index file.

\item \textbf{Stage}: Test the preconditioners.\\
The tests of the generated preconditioners must be enqueued manually by the user since no guarantee can be made, that storage quota will not be reached during the
generation phase. The \verb+--dependency=singleton+ option for SLURM has been used to ensure that the tests of generated preconditioners
are started only after all repetitions of the generation script have been run.
The testing stage produces, for each parameter set (experiment) and each repetition a unique text file containing the results of the 
execution of the chosen iterative method.
\end{enumerate}

Execution on the K80s differs in so far as the second stage is omitted and the first one is executed sequentially by a dedicated Python script
into which all the required parameters are hard-coded.

\subsection{Evaluation and expected result}

Evaluation of the numerical experiments is carried out by first consolidating the results into a single Pandas data frame.
This is done automatically by a preprocessing script which utilises the index of experiments generated in the first stage of the experiments.
The collected data is stored in CSV format.
It is then imported into a Jupyter notebook and further evaluation and visualization is performed in accordance with the requirements documented therein.

Raw results include plain-text output files from the SLURM scheduler and the code used to test the preconditioners.
Intermediate results are consolidated into CSV files and final results consist of a collection of plots showing the speed-up and execution time of 
different parameter configurations for different matrices. The images are stored in EPS format.

\subsection{Notes}
The MSPAI preconditioner may be obtained at https://www5.in.tum.de/wiki/index.php/MSPAI and is compiled with the provided Makefiles, which require 
the ATLAS library.

\end{document}